\documentclass[12pt]{amsart}
\newcommand{\deleted}[1]{}
\newcommand{\delete}[1]{}
\newcommand{\mynotes}[1]{}
\newcommand\notes[1]{}
\usepackage{txfonts}
\usepackage{amssymb}
\usepackage{eucal}
\usepackage{graphicx}
\usepackage{amsmath}
\usepackage{amscd}
\usepackage[all]{xy}
\usepackage[pdftex,  
            pdfstartview=FitH,
            CJKbookmarks=true,
            bookmarksnumbered=true,
            bookmarksopen=true,
            colorlinks, 
            pdfborder=001,   
            linkcolor=blue,
            anchorcolor=blue,
            citecolor=blue
            ]{hyperref}

\usepackage{amsthm}
\usepackage{mathrsfs}
\usepackage{enumerate}
\usepackage{amsfonts,latexsym}
\usepackage{xspace}
\usepackage{epsfig}
\usepackage{float}
\usepackage{color}
\usepackage{fancybox}
\usepackage{colordvi}
\usepackage{multicol}
\usepackage{tikz}
\usepackage{wasysym}
\usepackage{xspace}
\usepackage[active]{srcltx} 
\usepackage{mathtools} 
\usepackage{tabularx}
\usepackage{stfloats}

\newcommand\changed[1]{#1}
\changed{}

\newtheorem{theorem}{Theorem}[section]
\newtheorem{lemma}[theorem]{Lemma}

\newtheorem{coro}[theorem]{Corollary}

\newtheorem{prop}[theorem]{Proposition}
\theoremstyle{definition}

\newtheorem{remark}[theorem]{Remark}

\deleted{\newtheorem{theorem}{Theorem}[section]
\newtheorem{prop-def}{Proposition-Definition}[section]
\newtheorem{coro-def}{Corollary-Definition}[section]

}
\newcommand{\nc}{\newcommand}

\newcommand{\efootnote}[1]{}

\nc{\mlabel}[1]{\label{#1}}  
\nc{\mcite}[1]{\cite{#1}}  
\nc{\mref}[1]{\ref{#1}}  
\nc{\mbibitem}[1]{\bibitem{#1}} 
\renewcommand\geq{\geqslant}
\renewcommand\leq{\leqslant}

\renewcommand\bar[1]{\overline{#1}}
\newcommand{\lra}{\longrightarrow}

\newcommand{\cpx}[1]{#1^{\bullet}}
\newcommand{\stmodcat}[1]{#1\mbox{{\rm -{\underline{mod}}}}}
\newcommand{\Db}[1]{{\mathscr D}^b(#1)}
\newcommand{\Kb}[1]{{\mathscr K}^b(#1)}
\newcommand{\modcat}{\ensuremath{\mbox{{\rm -mod}}}}

\nc{\mrm}[1]{{\rm #1}}
\nc{\Hom}{\mrm{Hom}}
\nc{\Ext}{\mrm{Ext}}
\nc{\End}{\mrm{End}}
\nc{\rad}{\mrm{rad}}
\nc{\Aut}{\mrm{Aut}}
\nc{\X}{X^{\bullet}}
\nc{\Y}{Y^{\bullet}}
\nc{\Z}{Z^{\bullet}}
\nc{\T}{T^{\bullet}}
\nc{\Lv}{{\bf L}\nu}

\nc{\dcp}{the double centralizer property }
\nc{\ddcp}{the derived double centralizer property}
\nc{\add}{\mrm{add\,}}
\nc{\RHom}{\mrm{{\bf R}Hom}}
\nc{\OT}{\otimes^{\bf L}}
\nc{\thick}{\mrm{thick}}
\nc{\cone}{\mrm{cone}}
\nc{\proj}{\mrm{-proj}}
\nc{\Proj}{\mrm{Proj\,}}
\nc{\gl}{\mrm{gl}}
\nc{\im}{\mrm{im}}
\nc{\id}{\mrm{id}}
\nc{\redt}[1]{\textcolor{red}{#1}}
\nc{\jin}[1]{\textcolor{red}{Jing:#1}} 
\nc{\xing}[1]{\textcolor{purple}{Xing:#1}}

\begin{document}
\renewcommand{\thefootnote}{\alph{footnote}}
\setcounter{footnote}{-1} \footnote{ $^*$ Corresponding author.} 
\renewcommand{\thefootnote}{\alph{footnote}}
\setcounter{footnote}{-1}
\footnote{2020 Mathematics Subject Classification: Primary 18G80, 16E35, 16D50, Secondary 16G10, 18F30, 16P10.}
\renewcommand{\thefootnote}{\alph{footnote}}
\setcounter{footnote}{-1}
\footnote{Keywords: Derived equivalence, Grothendieck group, Nakayama permutation, self-injective algebra, symmetric algebra, weakly symmetric algebra.}

\title[Self-injective algebras under derived equivalences ]{Self-injective algebras under derived equivalences}

\author[Changchang Xi]{Changchang Xi}
\address{Changchang Xi: School of Mathematical Sciences, Capital Normal University, 100048 Beijing, \&
School of Mathematics and Statistics, Shaanxi Normal University, 710119 Xi'an, P. R. China}
\email{xicc@cnu.edu.cn}
\author[Jin Zhang]{Jin Zhang$^{*}$}
\address{Jin Zhang: School of Mathematical Sciences, Capital Normal University, 100048 Beijing, P. R. China}
\email{zj\_10@lzu.edu.cn}
\hyphenpenalty=8000
\begin{abstract}
The Nakayama permutations of two derived equivalent, self-injective Artin algebras are conjugate. A different but elementary approach is given to
showing that the weak symmetry and self-injectivity of finite-dimensional algebras over an \emph{arbitrary} field are preserved under derived equivalences.
\end{abstract}
\maketitle
\allowdisplaybreaks

\section{Introduction}

Self-injective algebras have played a very important role in various areas of mathematics and physics (see, for example, \cite{Ab}, \cite{kock}, and the references therein). In the representation theory of algebras, the famous, but still not yet solved  Auslander-Reiten conjecture on stable equivalences is reduced to self-injective Artin algebras (see \cite{Mar}). The conjecture states that two Artin algebras have same number of isomorphism classes of non-projective simple modules whenever they are stably equivalent. Moreover, self-injective Artin algebras are closed under stable equivalences (with a mild condition) (see \cite{Re}). Also, self-injective algebras over an algebraically closed field are closed under derived equivalences \cite{AR1}.  This result seems to be extended to algebras over any field in \cite{rr}, but we have difficulty to understand some arguments in its proof there. On the other hand, weakly symmetric algebras over an algebraically closed field are preserved under derived equivalences \cite{ad}. It is unknown whether this is true for algebras over an arbitrary field, while symmetric algebras over any field are closed under derived equivalences (see \cite{R}).

Given a self-injective Artin algebra $\Lambda$, the Nakayama functor of $\Lambda$ is a self-equivalence on the category of finitely generated projective $\Lambda$-modules. Hence it permutes the complete set of isomorphism classes of indecomposable projective $\Lambda$-modules. This permutation is called the Nakayama permutation, which is uniquely determined by $\Lambda$, up to conjugation.

This note continues the study of self-injective algebras under derived equivalences. First, we show that, if two self-injective Artin algebras are derived equivalent, then their Nakayama permutations are conjugate (see Theorem \ref{conj}). We then give an elementary approach to Rickard-Rouquier's result that self-injective algebras over an arbitrary field are closed under derived equivalences, and we further prove that derived equivalences also preserve weakly symmetric algebras over any field (see Corollary \ref{weasym}).

The strategy for proving the first result uses an idea from categorification, namely we first investigate derived Nakayama functors on the homotopy categories of finitely generated projective modules, and then pass to the Grothendieck groups of these homotopy categories. In this way,  the Nakayama permutations can be realized by derived Nakayama functors. The elementary proof of Rickard-Rouquier's result is based on studying relations between self-injective algebras and extensions of fields. Consequently, we prove the desired result for weakly symmetric algebras and self-injective algebras.

The paper is organized as follows: In Section \ref{sect2} we fix notation and recall basic facts on derived equivalences. In Section \ref{Groth} we study Grothendieck groups of triangulated categories. In Section \ref{Naka} we prove that the Nakayama permutations of derived equivalent, self-injective Artin algebras are conjugate. Also, we point out that Rickard's result on derived equivalences preserving symmetry for finite-dimensional algebras can be generalized to the one for Artin algebras (see Remark \ref{rmk4.3}).
In Section \ref{self-injective} we show that finite-dimensional, weakly symmetric algebras over an arbitrary field are closed under derived equivalences, and provide an elementary proof of Rickard-Rouquier's result on self-injective algebras under derived equivalences. Finally, we deduce a series of consequences of our main results.

\section{Preliminaries\label{sect2}}
In this section we fix notation and recall some definitions and results on derived equivalences.

Throughout the paper, all modules are assumed to be left modules. For a (unitary associative) ring $\Lambda$, we denote by $\Lambda$-mod the category of finitely generated $\Lambda$-modules, by $\Lambda$-proj the full subcategory of $\Lambda$-mod consisting of projective $\Lambda$-modules, and by $\Kb{\Lambda\mbox{\proj}}$  the bounded homotopy category of complexes over $\Lambda$-proj. As usual, we write $\Db{\Lambda}$ for the bounded derived category of $\Lambda\textmd{-mod}$.

Artin algebras $A$ and $B$ are \emph{derived equivalent} if $\Db{A}$ and $\Db{B}$ are equivalent as triangulated categories.
Derived equivalences can be described by tilting complexes \cite{R0}. We recall the descriptions just for Artin algebras below.

Let $A$ be an Artin algebra. A complex $\X$ in $\Kb{A\mbox{\proj}}$ is called a \emph{tilting complex} (see \cite{R0}) if $\Hom_{\Kb{A\mbox{\proj}}}(\X, \X[i])=0$ for $i\neq 0$ and $\X$ generates $\Kb{A\mbox{\proj}}$ as a triangulated category.

The description of derived equivalences by tilting complexes is given by the following theorem (see \cite{H, K1, R0, R}).

\begin{theorem} \label{derequ}
Suppose that $A$ and $B$ are Artin algebras over a commutative Artin ring $R$. Then the following are equivalent.

$(1)$ $A$ and $B$ are derived equivalent.

$(2)$ There exists a tilting complex $\T\in\Kb{A\emph{-proj}}$ such that
$B\simeq\End_{\Db{A}}(\T)^{op}$ as algebras.

$(3)$ There is a triangle equivalence from $\Kb{A\emph{-proj}}$ to $\Kb{B\emph{-proj}}$.
\end{theorem}

\medskip
An Artin algebra $A$ is said to be \emph{symmetric} if ${}_AA_A\simeq DA$ as $A$-$A$-bimodules, where $D$ is the usual duality of the Artin algebra $A$; \emph{weakly symmetric} if the injective hull and projective cover of every simple $A$-module are isomorphic; \emph{Frobenius} if ${}_AA\simeq DA$ as $A$-modules; and \emph{self-injective} if ${}_AA$ is injective. A basic self-injective algebra is a Frobenius algebra. By a \emph{basic} algebra we mean an Artin algebra $A$ such that $_AA$ is a direct sum of pairwise non-isomorphic indecomposable modules.

\medskip
Let $n$ be a positive integer. We denote by $\Sigma_n$ the symmetric group of permutations on $\{1,2, \cdots, n\}$. For an object $X$ in an additive category, $X^{\oplus n}$ stands for the direct sum of $n$ copies of $X$.

\section{Grothendieck groups of triangulated categories}\label{Groth}
In this section we study basic properties of the Grothendieck groups of triangulated categories, and their behaviors under triangle equivalences. We start with the following definition in \cite{Gro, H}.

Let $\mathcal{C}$ be a triangulated category with the shift functor $[1]$. Assume further that $\mathcal{C}$ is essentially small, that is, the isomorphism classes of objects of $\mathcal{C}$ form a set. For $X\in \mathcal{C}$, we denote by $[X]$ the isomorphism class containing $X$. Let $\widetilde{\mathcal{C}}$ be the set of the isomorphism classes $[X]$ of objects $X$ in $\mathcal{C}$. Let $\texttt{F}(\mathcal{C})$ be the free abelian group generated by all elements of $\widetilde{\mathcal{C}}$, and let $\texttt{F}_0(\mathcal{C})$ be the subgroup of $\texttt{F}(\mathcal{C})$ generated by $[X]-[Y]+[Z]$ for all triangles
$$X\lra Y\lra Z\lra X[1]$$
in $\mathcal{C}$. The \textit{Grothendieck group} $K_0(\mathcal{C})$ of $\mathcal{C}$ is defined to be the quotient group $\texttt{F}(\mathcal{C})/\texttt{F}_0(\mathcal{C})$. We write $\overline{[X]}$ for the coset of $[X]$ in $K_0(\mathcal{C})$.

We denote by
$$d: \widetilde{\mathcal{C}}\lra K_0(\mathcal{C})$$
the composition of the canonical maps $\widetilde{\mathcal{C}}\hookrightarrow \texttt{F}(\mathcal{C})\twoheadrightarrow K_0(\mathcal{C})$. Then $d([X])=\bar{[X]}$ for any object $X$ in $\mathcal{C}$.

For a triangle functor $F:\mathcal{C}\rightarrow \mathcal{D}$ of essentially small triangulated categories $\mathcal{C}$ and $\mathcal{D}$, one has naturally a map $\widetilde{F}: \widetilde{\mathcal{C}}\rightarrow \widetilde{\mathcal{D}}$ defined by $\widetilde{F}([X])=[F(X)]$ for  $[X]$ in $\widetilde{\mathcal{C}}$. Since the images of two isomorphic objects in $\mathcal{C}$ under $F$ are still isomorphic in $\mathcal{D}$, the map $\widetilde{F}$ is well defined.

\smallskip
If $A$ is a ring and $\mathcal{C}=\Kb{A\mbox{\proj}}$, we simply write $K_0(A)$ for $K_0(\mathcal{C})$.

\begin{lemma}\label{add} Let $\mathcal{C}$ be an essentially small triangulated category.

$(1)$ For objects $X$ and $Y$ in $\mathcal{C}$,
  $$d([X\oplus Y])=d([X])+d([Y]) \; \mbox{ and } \; d([X[i]])=(-1)^id([X]) \; \mbox{ for } \; i\in \mathbb Z.$$

$(2)$ The map $d$ is surjective.
\end{lemma}

{\it Proof.}
For objects $X$ and $Y$ in $\mathcal{C}$, there is a canonical triangle
$X \rightarrow X\oplus Y\rightarrow Y\stackrel{0} \rightarrow  X[1]$. Thus $d([X\oplus Y])=d([X])+d([Y])$. In particular, for the zero object $0$ of $\mathcal{C}$, there holds $d([0])=0$ in $K_0(\mathcal{C})$. The triangle
$X \rightarrow 0\rightarrow X[1]\stackrel{-\id_{X[1]}} \rightarrow  X[1]$ shows $d([X[1]])=-d([X])$. This implies $d([X[i]])=(-1)^id([X])$
for $i\in \mathbb Z$.

Let $\alpha\in K_0(\mathcal{C})$. Without loss of generality, we may assume that $\alpha$ is the coset of an element
$$r_1[X_1]+\cdots +r_m[X_m]+r_{m+1}[X_{m+1}]+\cdots +r_n[X_n]$$
in $\texttt{F}(\mathcal{C})$, where all $X_j$ are objects in $\mathcal{C}$, $r_j<0$ for $1\leq j\leq m$, and $r_j>0$ for $m+1\leq j\leq n$. Then
$$d([X_1[1]^{\oplus -r_{1}}\oplus \cdots \oplus X_m[1]^{\oplus -r_{m}}\oplus X_{m+1}^{\oplus r_{m+1}}\oplus\cdots \oplus X_n^{\oplus r_n}] )=\alpha.$$
So $d$ is surjective.
$\square$

\medskip
Now, assume that $A$ is a semiperfect ring, that is, every finitely generated left (or right) $A$-module has a projective cover, or equivalently, $A$ has a complete orthogonal set $\{e_1, \cdots, e_n\}$ of idempotents with each $e_iAe_i$ a local ring.

Let  $\X\in \Kb{A\mbox{\proj}}$ be of the form
$$ \X = \qquad \cdots \lra 0 \lra X^i\stackrel{d_X^{i}}\lra X^{i+1}\stackrel{d_X^{i+1}}\lra\cdots \lra X^{i+m}\lra 0 \lra \cdots$$
Denote by $_{\sigma{\leq t}}\X$ the brutal truncation of $\X$ at the degree $t$, that is, $(_{\sigma{\leq t}}\X)^j=X^j$ for $j\leq t$ and $0$ otherwise. Then there is a series of triangles:
$$
\Delta_{j}: \quad {_{\sigma{\leq j-1}}\X}[-1]\stackrel{f^{\bullet}_{j}}\longrightarrow X^{j}[-j]\lra {_{\sigma{\leq j}}\X}\lra {_{\sigma{\leq j-1}}\X}$$
for $ i+1\le j\le i+m$, where ${_{\sigma{\leq i+m}}\X}=\X$ and $f^{\bullet}_{j}$ is defined by $f^{j}_{j}=d_X^{j-1}$ and $f^{s}_{j}=0$ for $s\neq j$. By Lemma \ref{add}, $d([\X])=\sum_{i\leq j\leq i+m}(-1)^{j}d([X^j])$.

Since $A$ is a semiperfect ring, there are finitely many pairwise non-isomorphic, indecomposable projective $A$-modules. Let  $\{P_1, \cdots, P_n\}$ be a complete set of pairwise non-isomorphic, indecomposable projective $A$-modules. For $i\le j\le i+m$, we write $X^j\simeq \bigoplus_{1\leq s\leq n} P_s^{\oplus t_{j_s}}$ with $t_{j_s}\in \mathbb{N}$, and $\lambda_s :=\sum_{i\leq j\leq i+m}(-1)^jt_{j_s}$. Then $d([\X])=\sum_{1\leq s\leq n}\lambda_sd([P_s])$. As $d$ is surjective, we see that $K_0(A)$ is an abelian group generated by these $\bar{[P_s]}$, $1\le s\le n$.

Next, we define $\textbf{dim}(\X)=(\lambda_1, \cdots, \lambda_n)\in \mathbb{Z}^n$ for the complex $\X$. If $\Y$ is a complex in $ \Kb{A\mbox{\proj}}$ such that $\Y\simeq \X$ in $\Kb{A\mbox{\proj}}$, then $\textbf{dim}(\Y)=\textbf{dim}(\X)$. Moreover, for any morphism $\cpx{f}: \X\rightarrow \Z$ in $\Kb{A\mbox{\proj}}$, there holds $\textbf{dim}(\cone(\cpx{f}))=\textbf{dim}(\Z)-\textbf{dim}(\X)$, where $\cone(\cpx{f})$ stands for the mapping cone of $\cpx{f}$. Hence we get a homomorphism of abelian groups:
$$\textbf{dim}: K_0(A)\lra \mathbb{Z}^n,\quad \bar{[\X]}\mapsto \textbf{dim}(\X).$$This shows that the set $\{\textbf{dim}(\bar{[P_1]}), \textbf{dim}(\bar{[P_2]}), \cdots, \textbf{dim}(\bar{[P_n]})\}$ forms a basis of the free abelian group $\mathbb{Z}^n$, and therefore $K_0(A)$ is a free abelian group generated by $\bar{[P_1]},\bar{[P_2]}, \cdots, \bar{[P_n]}$.

\medskip
Not all Grothendieck groups of (essentially small) triangulated categories are free. For example, if $A$ is a finite-dimensional, self-injective algebra such that the Cartan matrix of $A$ has an elementary divisor different from $0$ and $1$, then the Grothendieck group of the stable module category $A\stmodcat$ (as a triangulated category) is not free. For more details, see \cite[Section 5.7.1]{Z}.

\begin{prop}\label{bar}
Suppose that $A$ and $B$ are semiperfect rings. If $F:\Kb{A\mbox{\proj}}\rightarrow \Kb{B\mbox{\proj}}$ is a triangle equivalence, then $F$ induces a group isomorphism $\bar{F}: K_0(A)\rightarrow K_0(B)$ such that the diagram (of maps) is commutative:
\[
\begin{array}{c}
(*)\qquad \xymatrix{
\widetilde{\Kb{A\mbox{\proj}}}\ar[r]^{\widetilde{F}}\ar[d]^{d}&\widetilde{\Kb{B\mbox{\proj}}}
\ar[d]^{d}\\
K_0(A)\ar[r]^{\bar{F}}&K_0(B)
}
\end{array}
\]where $\widetilde{\Kb{A\mbox{\proj}}}$ stands for the set of the isomorphism classes of objects in $\Kb{A\mbox{\proj}}.$
\end{prop}

{\it Proof.}
Let $\Y_s:=F(P_s)$ in $\Kb{B\mbox{\proj}}$ for $1\leq s\leq n$. We define a group homomorphism $$\bar{F}: K_0(A)\lra K_0(B),\quad \bar{[P_s]}\mapsto d([\Y_s])\; \mbox{  for } 1\le s\le n.$$Let $\X$ be a complex in $\Kb{A\mbox{\proj}}$ of the above form with $X^j  := \bigoplus_{1\leq s\leq n} P_s^{\oplus t_{j_s}}$ and $\lambda_s:=\sum_j(-1)^jt_{j_s}$. Then $d([\X])=\sum_{1\leq s\leq n}\lambda_sd([P_s])$ and $$\bar{F}d([\X])=\sum_{1\leq s\leq n}\lambda_s\bar{F}d([P_s])=\sum_{1\leq s\leq n}\lambda_s\bar{F}(\bar{[P_s]})=\sum_{1\leq s\leq n}\lambda_sd([\Y_s]).$$

Applying $F$ to the triangle $\Delta_j$, we have a triangle
$$
F(\Delta_{j}): \quad F({_{\sigma{\leq j-1}}\X})[-1]\stackrel{F(f^{\bullet}_{j})}\lra F(X^{j})[-j]\lra F({_{\sigma{\leq j}}\X})\lra F({_{\sigma{\leq j-1}}\X}),\;\; i+1\le j\le i+m.$$
Since $F({_{\sigma{\leq i+m}}\X})=F(\X)$, it follows from Lemma \ref{add} that $$d([F(\X)])=\sum_{i\leq j\leq i+m}(-1)^{j}d([F(X^j)]).$$ Then
$$d\widetilde{F}([\X])=d([F(\X)])=\sum_{1\leq s\leq n}\lambda_sd([F(P_s)])=\sum_{1\leq s\leq n}\lambda_sd([\Y_s])=\bar{F}d([\X]).$$
Hence the above square $(*)$ is commutative.

It remains to show that $\bar{F}$ is bijective. In fact, we consider a quasi-inverse $F^{-1}$ of $F$. In this case, we have the group homomorphism $\bar{F^{-1}}: K_0(B)\to K_0(A)$ induced by $F^{-1}$. Then
$$\bar{F^{-1}}\,\bar{F}(\bar{[P_s]})= \bar{F^{-1}}(d([\Y_s])) =d\widetilde{F^{-1}}([\Y_s])=d([F^{-1}(\Y_s)])=d([P_s])= \bar{[P_s]}$$for $1\leq s\leq n$. Hence $\bar{F^{-1}}\,\bar{F}=\id_{K_0(A)}$. Similarly, $\bar{F}\,\bar{F^{-1}}=\id_{K_0(B)}$. So $\bar{F}$ is bijective. $\square$

\section{Nakayama permutations of self-injective algebras}\label{Naka}

In this section we will prove that Nakayama permutations of derived equivalent, self-injective Artin algebras are conjugate.

Let $A$ be an Artin algebra over a commutative Artin ring $R$. The Nakayama functor $\nu_A: A\textmd{-mod}\rightarrow A\textmd{-mod}$ is defined by
$\nu_A:=(DA)\otimes_A-,$
where $D$ is the usual duality of an Artin algebra. Clearly, $\nu_A$ induces a left derived functor $$\Lv_A:\Db{A}\lra \Db{A},$$which restricts to a triangle equivalence
$$\Lv_A: \Kb{A\mbox{\proj}}\lra  \Kb{A\mbox{-inj}}, $$where $A\mbox{-inj}$ denotes the category of finitely generated injective $A$-modules.

\medskip
Now, assume that $A$ is a self-injective Artin algebra. Since $_A(DA)$ and $(DA)_A$ are projective generators, $\nu_A$ is a self-equivalence on $A$-mod, and restricts to a self-equivalence on $A\mbox{\proj}$. Let $\{P_1, \cdots, P_n\}$ be a complete set of pairwise non-isomorphic, indecomposable projective $A$-modules. Then $\nu_A$ induces a permutation on $\{P_1, \cdots, P_n\}$, called the Nakayama permutation of $A$. Precisely, the \textit{Nakayama permutation} $\sigma_A$ is defined on $\{1, \cdots, n\}$ by
$$\nu_A(P_i)\simeq P_{\sigma_A(i)}$$
for $i\in \{1, \cdots, n\}$. Clearly, up to conjugation, the Nakayama permutation $\sigma_A$ of $A$ is uniquely determined by $\{P_1, \cdots, P_n\}$.

Let $B$ be another self-injective Artin algebra over $R$, and let $\{Q_1, \cdots, Q_m\}$ be a complete set of pairwise non-isomorphic, indecomposable projective $B$-modules. Assume that $A$ and $B$ are derived equivalent. Then $m=n$ and $\sigma_B$ is again a permutation of $\{1, \cdots, n\}$. Our first main result reveals a precise relation between $\sigma_A$ and $\sigma_B$.

\begin{theorem}\label{conj}
If $A$ and $B$ are derived equivalent, self-injective Artin algebras, then $\sigma_A$ and $\sigma_B$ are conjugate.
\end{theorem}

To prove  Theorem \ref{conj}, we first show a technical lemma on the left derived functors of Nakayama functors.

Let $A$ be a self-injective Artin algebra. Then both $_A(DA)$ and $(DA)_A$ are projective. By definition, the left derived functor of the Nakayama functor $\nu_A$ is given explicitly as follows:
$$\Lv_A: \Kb{A\mbox{\proj}}\lra  \Kb{A\mbox{\proj}}, \;\; \X=(X^i, d^i_X)\mapsto \big(\nu_A(X^i), \nu_A(d^i_X)\big).$$
As $\nu_A$ is a self-equivalence of $A\mbox{\proj}$ (or $A$-mod), we see that $\Lv_A$ is a triangle self-equivalence of $\Kb{A\mbox{\proj}}$. Clearly, $\Lv_A(P)\simeq \nu_A(P)$ for $P\in A\mbox{\proj}$.

\begin{lemma}\label{Lv}
 Suppose that $A$ and $B$ are Artin algebras. If $F: \Db{A}\rightarrow  \Db{B}$ is a triangle equivalence, then for any $\X$ in $\Kb{A\mbox{\proj}}$, $F\Lv_A(\X)\simeq \Lv_BF(\X)$ in $\Db{B}$ which is natural in $\X$. In particular, if $A$ and $B$ are self-injective Artin algebras and $F: \Kb{A\mbox{\proj}}\rightarrow  \Kb{B\mbox{\proj}}$ is a triangle equivalence, then there is a natural isomorphism $F\Lv_A\simeq \Lv_BF: \Kb{A\mbox{\proj}}\rightarrow  \Kb{B\mbox{\proj}}$.
\end{lemma}

{\it Proof.}
For $\X$ in $\Kb{A\mbox{\proj}}$ and $\Y$ in $\Db{A}$, we may consider $\Hom_{\Db{A}}(\X, \Y)$ and \\$\Hom_{\Db{A}}(\Y, \Lv_A(\X))$ as the degree-zero homologies of the total complexes of the double complexes $\Hom_A^{\bullet\bullet}(\X, \Y)$ and $\Hom_{A}^{\bullet\bullet}(\Y, \Lv_A(\X))$, respectively. It is well known that, for any $X$ in $A\mbox{\proj}$ and $Y$ in $A\mbox{\textmd{-mod}}$, $D\Hom_A(X, Y)\simeq \Hom_A(Y, \nu_A(X))$ which is natural in $X$ and $Y$. Thus $D\Hom_A^{\bullet\bullet}(\X, \Y)\simeq\Hom_{A}^{\bullet\bullet}(\Y, \Lv_A(\X))$ naturally as double complexes for $\X\in \Kb{A\mbox{\proj}}$ and $\Y\in \Db{A}$. Taking homology in degree zero, we obtain
$$(1)\quad D\Hom_{\Db{A}}(\X, \Y)\simeq \Hom_{\Db{A}}(\Y, \Lv_A(\X))$$
which is natural in $\X\in \Kb{A\mbox{\proj}}$ and $\Y\in \Db{A}.$
On the other hand, as $F$ is an equivalence, there are natural isomorphisms:
$$(2) \quad D\Hom_{\Db{A}}(\X, \Y)\simeq D\Hom_{\Db{B}}(F(\X), F(\Y))$$ and $$(3) \quad\Hom_{\Db{A}}(\Y, \Lv_A(\X))\simeq \Hom_{\Db{B}}(F(\Y), F(\Lv_A(\X))).$$ Using the $B$-module version of $(1)$, we have the natural isomorphism: $$(4) \quad D\Hom_{\Db{B}}(F(\X), F(\Y))\simeq \Hom_{\Db{B}}(F(\Y), \Lv_B(F(\X))).$$ Thus it follows from $(3), (1), (2)$ and $(4)$ that
$$\Hom_{\Db{B}}(F(\Y), F(\Lv_A(\X)))\simeq\Hom_{\Db{B}}(F(\Y), \Lv_B(F(\X)))$$ which is natural in $\X\in \Kb{A\mbox{\proj}}$ and $\Y\in \Db{A}$. As $F$ is an equivalence, we obtain $F\Lv_A(\X)\simeq \Lv_BF(\X)$ in $\Db{B}$ which is natural in $\X\in \Kb{A\mbox{\proj}}$. When $A$ and $B$ are self-injective and when $F$ restricts to a triangle equivalence $\Kb{A\mbox{\proj}}\to  \Kb{B\mbox{\proj}}$, we have the last statement of Lemma \ref{Lv}. $\square$

\medskip
\begin{remark} \label{rmk4.3} Lemma \ref{Lv} can be applied to generalize \cite[Corollary 5.3]{R} for finite-dimensional algebras over a field to the one for Artin algebras, namely an Artin algebra $B$ derived equivalent to a symmetric Artin algebra $A$ is itself symmetric. Indeed, let $F: \Db{A}\rightarrow \Db{B}$ be a triangle equivalence. Since $A$ is symmetric, $DA\simeq A$ as $A$-$A$-bimodules, and therefore $\Lv_A\simeq \id$ naturally on $\Db{A}$. By Lemma \ref{Lv}, $\Lv_BF(\X)\simeq F\Lv_A(\X)\simeq F(\X)$ in $\Db{B}$  naturally for $\X$ in $\Kb{A\mbox{\proj}}$. As $F$ is an equivalence, $\Lv_B\simeq \id$ naturally on $\Kb{B\mbox{\proj}}$. Hence $DB\simeq B$ as $B$-modules.
If we apply the natural isomorphism $\Lv_B\simeq \id$ to morphisms $B\rightarrow B$ in $\Kb{B\mbox{\proj}}$ given by right multiplication of elements in $B$, then the isomorphism $DB\simeq B$ is actually an isomorphism of $B$-$B$-bimodules, and therefore $B$ is a symmetric algebra.
\end{remark}
\medskip
Given $\sigma\in \Sigma_n$, we may write $\sigma=\sigma_1\sigma_2\cdots\sigma_s$ with $\sigma_i$ a cyclic permutation of length $\lambda_i\ge 1$, such that the contents of these $\sigma_i$ are pairwise disjoint. In this case, we may assume that $\lambda_1\geq \lambda_2\geq \cdots \geq \lambda_s$. Then $\lambda:=(\lambda_1, \cdots, \lambda_s)$ is a partition of $n$, called the \emph{cycle type} of $\sigma$.
It is well known that two permutations in $\Sigma_n$ are conjugate if and only if they have the same cycle type.

For a unitary ring $A$, we denote by $M_n(A)$ the full $n\times n$ matrix ring over $A$. If $\sigma\in \Sigma_n$, then \emph{the permutation matrix} $c_{\sigma}$ of $\sigma$ over $\mathbb{C}$ is the $n\times n$ matrix with $1$ in the $(i,\sigma(i))$-entry for $1\le i\le n$ and with $0$ for all other entries.

The following result seems to be known. For the convenience of the reader, we provide a proof.

\begin{lemma}\label{per}
Let $\sigma_1$ and $\sigma_2$ be permutations in $\Sigma_n$. Then $\sigma_1$ and $\sigma_2$ are conjugate in $\Sigma_n$ if and only if $c_{\sigma_1}$ and $c_{\sigma_2}$ are similar in $M_n(\mathbb{C})$.
\end{lemma}

{\it Proof.} Clearly, $c_{\sigma_1}c_{\sigma_2}=c_{\sigma_1\sigma_2}$ and $c_{\sigma}^{-1}=c_{\sigma^{-1}}$ in $M_n(\mathbb{C})$. Thus, if $\sigma_1$ and $\sigma_2$ are conjugate in $\Sigma_n$, then $c_{\sigma_1}$ and $c_{\sigma_2}$ are similar in $M_n(\mathbb{C})$.
Here, $\mathbb{C}$ can be replaced by any field.

Conversely, suppose that $c_{\sigma_1}$ and $c_{\sigma_2}$ are similar in $M_n(\mathbb{C})$. Let $\lambda=(\lambda_1, \cdots, \lambda_u)$ and $\mu=(\mu_1, \cdots, \mu_v)$ be the cycle types of $\sigma_1$ and $\sigma_2$, respectively. Since the similarity of matrices and conjugation of permutations are equivalence relations and since conjugate permutations have similar permutation matrices, we may assume that
$$\sigma_1=(1\, 2\, \cdots\, \lambda_1)(\lambda_1+1\,\lambda_1+2\,\cdots\, \lambda_1+\lambda_2)\cdots (\sum_{1\leq i\leq u-1}\lambda_i+1\,\sum_{1\leq i\leq u-1}\lambda_i+2\,\cdots\, n),$$
$$\sigma_2=(1\, 2\, \cdots\, \mu_1)(\mu_1+1\,\mu_1+2\,\cdots\, \mu_1+\mu_2)\cdots (\sum_{1\leq i\leq v-1}\mu_i+1\,\sum_{1\leq i\leq v-1}\mu_i+2\,\cdots\, n)$$
where the $i$-tuple $(a_1\, \cdots\, a_i)$ means the cyclic permutation on the set $\{a_1, \cdots, a_i\}$.
By computations, the characteristic polynomials of $c_{\sigma_1}$ and $c_{\sigma_2}$ are
 $$\Phi_1(x)=(x^{\lambda_1}-1)(x^{\lambda_2}-1)\cdots (x^{\lambda_u}-1)\in \mathbb{C}[x]\; \mbox{ and } \; \Phi_2(x)=(x^{\mu_1}-1)(x^{\mu_2}-1)\cdots (x^{\mu_v}-1)\in \mathbb{C}[x],$$
respectively. Since $c_{\sigma_1}$ and $c_{\sigma_2}$ are similar in $M_n(\mathbb{C})$, we have $\Phi_1(x)=\Phi_2(x)$, that is, they have the same eigenvalues with the same multiplicities. We show $\lambda_1=\mu_1$. This follows from the $3$ facts:

(i) All $\lambda_1$-th roots of unity are eigenvalues of $c_{\sigma_1}$, while all $\mu_1$-th  roots of unity are eigenvalues of $c_{\sigma_2}$.

(ii) There exists a $q$-th root of unity different from any $w$-th root of unity if $q>w$, and

(iii) $\lambda_1$ and $\mu_1$ are maximal in $\lambda$ and $\mu$, respectively.

By repeating this process, we finally get $u=v$ and $\lambda_i=\mu_i$ for $1\le i\le u$. Hence $\sigma_1$ and $\sigma_2$ have the same cycle type, and therefore are conjugate in $\Sigma_n$.
$\square$

\medskip
{\bf Proof of Theorem \ref{conj}:}
The functors $\Lv_A$ and $\Lv_B$ are triangle self-equivalences of $\Kb{A\mbox{\proj}}$ and $\Kb{B\mbox{\proj}}$, respectively.
Since $A$ and $B$ are derived equivalent, there is a triangle equivalence $F:\Kb{A\mbox{\proj}}\rightarrow \Kb{B\mbox{\proj}}$ by Theorem \ref{derequ}. Thus there is the following diagram:
\[\xymatrix@C=0.7em@R=1em{
&&\widetilde{\Kb{A\mbox{\proj}}}\ar[rr]^{\widetilde{F}}\ar@{.>}[dd]_>>>>>>>>>>>>{d} &                         &\widetilde{\Kb{B\mbox{\proj}}}\ar[dd]^{d}                     \\
& \widetilde{\Kb{A\mbox{\proj}}}\ar[dd]^<<<<<<{d} \ar[rr]^>>>>>>>>>>{\quad \; \widetilde{F}}\ar[ur]^{\widetilde{\Lv_A}}   &   &\widetilde{\Kb{B\mbox{\proj}}}\ar[dd]^>>>>>>>>>>>{d} \ar[ur]^{\widetilde{\Lv_B}}   &       \\
&&K_0(A) \ar@{.>}[rr]^<<<<<<<<<{\bar{F}} &                         &K_0(B)\\
&K_0(A)\ar@{.>}[ur]^>>>>>>>>{\bar{\Lv_A}}\ar[rr]_{\bar{F}}    &    &K_0(B)\ar[ur]_>>>>>>>>{\bar{\Lv_B}} &\\
}
\]
where the vertical squares in the diagram are commutative by Proposition \ref{bar}, and the top square is commutative by Lemma \ref{Lv}. We shall show that the bottom square of homomorphisms of abelian groups is also commutative, that is, $\bar{F}\,\bar{\Lv_A}=\bar{\Lv_B}\,\bar{F}$.

Indeed, take $\alpha\in K_0(A)$. By Lemma \ref{add}(2), there is a complex $\X$ in $\Kb{A\mbox{\proj}}$ such that $d([\X])=\alpha$. Then
\begin{equation*}
 \begin{aligned}
\bar{F}\,\bar{\Lv_A}(\alpha)&=\bar{F}\,\bar{\Lv_A}(d([\X]))=\bar{F}d\widetilde{\Lv_A}([\X])=d\widetilde{F}\widetilde{\Lv_A}([\X])\\
&=d\widetilde{\Lv_B}\widetilde{F}([\X])=\bar{\Lv_B}d\widetilde{F}([\X])=\bar{\Lv_B}\,\bar{F}d([\X])=\bar{\Lv_B}\,\bar{F}(\alpha).
\end{aligned}
\end{equation*}
Hence the bottom square of the diagram is commutative.

Now, consider the Nakayama permutations $\sigma_A$ and $\sigma_B$ as elements in $\Sigma_n$, which are defined by $\nu_A(P_i)\simeq P_{\sigma_A(i)}$  and $\nu_B(Q_i)\simeq Q_{\sigma_B(i)}$ for $1\le i\le n$. Let $c_{\sigma_A}$ and $c_{\sigma_B}$ be the permutation matrices of $\sigma_A$ and $\sigma_B$, respectively. The Grothendieck groups $K_0(A)$ and $K_0(B)$ are free abelian groups generated by these $\bar{[P_i]}$ and these $\bar{[Q_i]}$, respectively. Moreover, by Proposition \ref{bar},
$$\bar{\Lv_A}(\bar{[P_i]})=\bar{\Lv_A}(d([P_i]))=d\widetilde{\Lv_A}([P_i])=d([\nu_A(P_i)]) = d([P_{\sigma_A(i)}])=\bar{[P_{\sigma_A(i)}]}.$$
Hence, with respect to the basis $\{\bar{[P_1]}, \cdots, \bar{[P_n]}\}$, the group homomorphism $\bar{\Lv_A}$ has the corresponding matrix $c_{\sigma_A}$. Similarly, with respect to the basis $\{\bar{[Q_1]}, \cdots, \bar{[Q_n]}\}$, the group homomorphism $\bar{\Lv_B}$ has the corresponding matrix $c_{\sigma_B}$. Since $\bar{F}$ is a group isomorphism by Proposition \ref{bar}, it corresponds to an invertible matrix $c\in M_n(\mathbb{C})$ with respect to the basis $\{\bar{[P_1]}, \cdots, \bar{[P_n]}\}$ of $K_0(A)$ and the basis $\{\bar{[Q_1]}, \cdots, \bar{[Q_n]}\}$ of $K_0(B)$. Due to $\bar{F}\,\bar{\Lv_A}=\bar{\Lv_B}\,\bar{F}$, there holds $cc_{\sigma_A}=c_{\sigma_B}c$. This means that $c_{\sigma_A}$ and $c_{\sigma_B}$ are similar in $M_n(\mathbb{C})$. By Lemma \ref{per},  $\sigma_A$ and $\sigma_B$ are conjugate in $\Sigma_n$.
$\square$

\section{Self-injective and Weakly symmetric algebras over a field are closed under derived equivalences }\label{self-injective}

Al-Nofayee and Rickard \cite{AR1} proved that, if $A$ and $B$ are derived equivalent, finite-dimensional algebras over an algebraically closed field and if $A$ is self-injective, then $B$ is self-injective. This result seems then to be extended to finite-dimensional algebras over an \emph{arbitrary} field by Rickard and Rouquier in \cite[Corollary 3.12]{rr}, but we have difficulty to understand an argument in the proof there, see the words just above \cite[Corollary 3.12]{rr}: ``Assume now $H^{<0}(B) =0$. Then, viewed as an object of $D^b(B), \nu(P_{\mathcal{S}}(S))$ is concentrated in degree $0$".

In this section we give a different, but very elementary approach to Rickard-Rouquier's result, and we show that a finite-dimensional algebra over an arbitrary field derived equivalent to a weakly symmetric algebra is itself weakly symmetric. This is known for weakly symmetric algebras over an algebraically closed field in \cite[Proposition 3.1]{ad}.

An Artin ring $R$ is called a \emph{Frobenius ring} if $_RR$ is injective and the socle of $_RR$ is isomorphic to the top of $_RR$.
\begin{lemma}\label{pro-inj}
Let $\Lambda$ be an Artin algebra over a Frobenius and commutative Artin ring $R$, and let $E$ be a commutative $R$-algebra such that $_RE$ is a free $R$-module and $_EE$ is an injective $E$-module. Assume that $M\in {\Lambda}\modcat$ is a projective $R$-module. Then $_{\Lambda}M$ is injective if and only if so is the $\Lambda\otimes_RE$-module $M\otimes_R E$.
\end{lemma}

{\it Proof.} Let ${}-^*=\Hom_R(-,R): \Lambda\textmd{-mod}\rightarrow \Lambda^{op}\textmd{-mod}$. Since $R$ is a Frobenius ring, $-^*$ is a duality by \cite[Theorem 3.3]{AR}. With $M$ also  $M^*$ is a finitely generated projective $R$-module, therefore $\Hom_R(M^*,R)\otimes_RX\simeq \Hom_R(M^*,X)$ as $\Lambda$-$\Gamma$-bimodules for any $R$-$\Gamma$-bimodule $X$ with $\Gamma$ a ring. Hence there are natural isomorphisms of functors:
$$\begin{array}{rl}
\Hom_{\Lambda\otimes_RE}(-, M\otimes_RE) &\simeq \Hom_{\Lambda\otimes_RE}(-, (M^{**})\otimes_RE)\\
                                         & \simeq \Hom_{\Lambda\otimes_RE}\big(-, \Hom_R(M^*,E)\big)\\
                                         & \simeq \Hom_E\big(M^*\otimes_{\Lambda}(-)_E, E\big) \qquad (\mbox{by adjoint isomorphism}) \\
                                         & = \Hom_E(-,E)\circ (M^*\otimes_{\Lambda}-)\\
\end{array}$$
Thus if $_{\Lambda}M$ is injective, then it follows from the duality $-^*$ that $M^*$ is a projective right $\Lambda$-module, and therefore $\Hom_{\Lambda\otimes_RE}(-, M\otimes_RE)$ is a composition of two exact functors. Thus  $\Hom_{\Lambda\otimes_RE}(-, M\otimes_RE)$ is itself an exact functor and $M\otimes_RE$ is an injective $\Lambda\otimes_RE$-module.

Conversely, suppose that $_{\Lambda\otimes_RE}(M\otimes_RE)$ is injective. By \cite[Corollary IX.2.4a]{CE}, the $\Lambda$-module $M\otimes_RE$ is injective. Assume that $\{x_i\mid i\in I\}$ is an $R$-basis of $E$ for some indexing set $I$. We take a fixed element $0\in I$. Then $E=x_0 R\oplus \bigoplus_{i\in I\setminus 0}x_iR$ and
$$_{\Lambda}M\otimes_RE\simeq M\otimes_R \big(\bigoplus_{i\in I}x_iR\big)\simeq M\otimes_R (x_0R)\oplus M\otimes_R \big(\bigoplus_{i\in I\setminus 0} x_iR\big)\simeq M\oplus M\otimes_R \big(\bigoplus_{i\in I\setminus 0}x_iR\big)$$as $\Lambda$-modules. Hence $_{\Lambda}M$ is injective.
$\square$

\medskip
The following is an immediate consequence of Lemma \ref{pro-inj}.
\begin{coro}\label{pro-inj}
Let $\Lambda$ be a finite-dimensional algebra over a field $k$, and let $E/k$ be an extension of fields. Then $\Lambda$ is self-injective if and only if so is the tensor product $\Lambda\otimes_k E$ of the $k$-algebras $\Lambda$ and $E$.
\end{coro}

The following result is observed by  Rickard and Rouquier in \cite{rr}.

\begin{coro}\label{eq} \cite[Corollary 3.12]{rr}
Suppose that $A$ and $B$ are finite-dimensional algebras over an arbitrary field such that they are derived equivalent. If $A$ is self-injective, then so is $B$.
\end{coro}

{\it Proof.} Assume that $A$ and $B$ are finite-dimensional algebras over a field $k$. Let $\bar{k}$ be an algebraically closed field of $k$. Since $A$ and $B$ are derived equivalent, $A \otimes_k \bar{k}$ and $B \otimes_k \bar{k}$ are derived equivalent by \cite[Theorem 2.1]{R}. Suppose that $A$ is self-injective. By Corollary \ref{pro-inj}, $A \otimes_k \bar{k}$ is self-injective. It is easy to see that $A \otimes_k \bar{k}$ and $B \otimes_k \bar{k}$ are finite-dimensional algebras over $\bar{k}$ because $\dim_{\bar{k}}(A \otimes_k \bar{k})=\dim_{k}(A)$. Now, the $\bar{k}$-algebra $B \otimes_k \bar{k}$ is self-injective by \cite[Theorem 2.1]{AR1} which states that finite-dimensional self-injective algebras over an algebraically closed field are preserved under derived equivalences. It then follows from Corollary \ref{pro-inj} that $B$ is a self-injective algebra.
$\square$

\medskip
Derived equivalences preserve finite-dimensional symmetric algebras over an arbitrary field \cite[Corollary 5.3]{R}. We point out that this is true also for weakly symmetric algebras over an arbitrary field, and refer to \cite{ad} for weakly symmetric algebras over an algebraically closed field.

\begin{coro}\label{weasym}
 Suppose that $A$ and $B$ are finite-dimensional algebras over an arbitrary field such that they are derived equivalent. If $A$ is weakly symmetric, then so is $B$.
\end{coro}

{\it Proof.} A finite-dimensional, self-injective algebra $\Lambda$ is weakly symmetric if and only if the Nakayama permutation of $\Lambda$ is the identity map. This follows from the definition of the Nakayama functor $\nu_{\Lambda}$.

Suppose that $A$ and $B$ are derived equivalent. Further, assume that $A$ is weakly symmetric. Then $A$ is self-injective. By Corollary \ref{eq}, $B$ is also self-injective. By assumption, $A$ is weakly symmetric, that is, the Nakayama permutation of $A$ is the identity map. So the Nakayama permutation of $B$ is also the identity map by Theorem \ref{conj}. Hence $B$ is weakly symmetric.
$\square$

\begin{coro}\label{tensor}
Suppose that $A,B, \Lambda$ and $\Gamma$ are finite-dimensional algebras over a field $k$. Assume that $A$ and $\Lambda$ are derived equivalent and that $B$ and $\Gamma$ are derived equivalent.
\begin{enumerate}
\item If both $A$ and $B$ are symmetric (or self-injective), then so is the tensor product algebra $\Lambda\otimes_k\Gamma$.
\item Assume that $k$ is an algebraically closed field. If both $A$ and $B$ are weakly symmetric, then so is the tensor product algebra  $\Lambda\otimes_k\Gamma$.
\end{enumerate}
\end{coro}

{\it Proof.} It is known that if $A$ and $B$ are symmetric (or self-injective) $k$-algebras over an arbitrary field $k$, then so is the tensor product $A\otimes_kB$. Furthermore, if $A$ and $B$ are weakly symmetric algebras over an algebraically closed  $k$, then so is the tensor product $A\otimes_kB$. This can be seen by the two facts: (a) Over an algebraically closed field $k$, the indecomposable projective $(A\otimes_kB)$-modules are of the form $P\otimes_kQ$ with $P$ and $Q$ indecomposable modules over $A$ and $B$, respectively. (b) Over an arbitrary field $k$, there holds $\nu_{A\otimes_kB}(P\otimes_kQ)\simeq \nu_A(P)\otimes_k\nu_B(Q)$ as $(A\otimes_kB)$-modules for $P$ in $A\mbox{\proj}$ and $Q$ in $B\mbox{\proj}$. As derived equivalences are preserved under taking tensor products (see \cite{R}), we see that $A\otimes_kB$ and $\Lambda\otimes_k\Gamma$ are derived equivalent. Now, Corollary \ref{tensor} follows from Remark \ref{rmk4.3} (or \cite[Corollary 5.3]{R}) and Corollaries \ref{eq}-\ref{weasym}. $\square$

\medskip
Let $\mathcal{C}$ be an additive category, $\mathcal{D}$ a full
subcategory of $\mathcal C$, and $X$ an object in $\mathcal{C}$. A morphism $f: D\to X$ in $\mathcal{C}$ is called a right $\mathcal{D}$-approximation of $X$ if $D\in \mathcal{D}$ and the induced map
$\Hom_{\mathcal{C}}(D',f): \Hom_{\mathcal{C}}(D', D)\to \Hom_{\mathcal{C}}(D',X)$ is surjective for every object $D'\in\mathcal{D}$. Dually, there is defined the left $\mathcal{D}$-approximation of $X$.

A sequence $$X\stackrel{f}{\longrightarrow}M\stackrel{g}{\longrightarrow}Y$$ in $\mathcal C$
is called a $\mathcal D$-\emph{split sequence} if $M\in {\mathcal D}$, $f$ is both a kernel of $g$ and a left $\mathcal D$-approximation of $X$, and  $g$ is both a cokernel of $f$ and a right $\mathcal D$-approximation of $Y$.

For an object $M$ in $\mathcal{C}$, $\add(M)$ stands for the full subcategory of $\mathcal{C}$ consisting of all objects isomorphic to direct summands of direct sums of finitely many copies of $M$.

\smallskip
As a consequence of Corollary \ref{eq} and Corollary \ref{weasym} together with \cite[Theorem 3.5]{hx2011} and \cite[Corollary 5.3]{R}, we get the following.

\begin{coro}\label{d-splseq} Let $M$ be an object of an additive $k$-category $\mathcal{C}$ with $k$ a field.
If $X\rightarrow M'\rightarrow Y$ is an $\add(M)$-split sequence in $\mathcal{C}$, then $\End_{\mathcal{C}}(X\oplus M)$ is a self-injective (symmetric, weakly symmetric) algebra if and only if so is $\End_{\mathcal{C}}(Y\oplus M)$.
\end{coro}

Let $A$ be an Artin algebra. A complex $\cpx{T}$ in $\Kb{A\mbox{\proj}}$ is called a \emph{basic complex} if it is a direct sum of pairwise non-isomorphic, indecomposable complexes in $\Kb{A\mbox{\proj}}$.
A complex $\X$ in $\Kb{A\mbox{\proj}}$ is said to be \emph{radical} if all differentials of $\X$ are radical homomorphisms.

\begin{coro}\label{basic}
If $A$ is a finite-dimensional, self-injective algebra over a field, then, for any basic tilting complex $\X$, $\Lv_A(\X)\simeq \X$ in $\Kb{A\mbox{\proj}}$.
\end{coro}

{\it Proof.} Let $\X$ be a basic tilting complex and $B:=\End_{\Kb{A\mbox{\proj}}}(\X)^{op}$. Then $B$ is a basic self-injective algebra by Corollary \ref{eq}, and therefore $B$ is a Frobenius algebra. By definition, $_B(DB)\simeq {}_BB$ as $B$-modules. Moreover, there is a triangle equivalence $F: \Kb{B\mbox{\proj}}\rightarrow \Kb{A\mbox{\proj}}$ such that $F(_BB)=\X$ (see \cite{R0}). Then it follows from Lemma \ref{Lv} that
$$\Lv_A(\X)=\Lv_A(F(_BB))\simeq F\Lv_B(_BB)=F(_B(DB))\simeq F(_BB)=\X$$in $\Kb{A\mbox{\proj}}$.
$\square$

\begin{coro}
Suppose that $A$ and $B$ are finite-dimensional algebras over a field such that they are derived equivalent. If $A$ is self-injective and its Nakayama permutation $\sigma_A$ is transitive (that is, $\sigma_A$ has only one orbit), then $A$ and $B$ are Morita equivalent.
\end{coro}

{\it Proof.}
Without loss of generality, we assume that $B$ is basic. Then there is a basic tilting complex $\X$ in $\Kb{A\mbox{\proj}}$ such that $B\simeq \End_{\Kb{A\mbox{\proj}}}(\X)^{op}$ by Theorem \ref{derequ}. Further, by \cite[(a), p.112]{HX}, we may assume that the complex $\X$ is radical. Now it suffices to prove that $\X$ is concentrated in a single degree because this will imply that $\X$ is a projective generator, and therefore $A$ and $B$ are Morita equivalent.

Indeed, assume that $\X$ is of the form (up to shift)
$$ \X = \qquad \cdots \lra 0 \lra X^0\stackrel{d_X^{0}}\lra X^{1}\stackrel{d_X^{1}}\lra\cdots \lra X^{m}\lra 0 \lra \cdots$$
with $X^0\ne 0\ne X^{m}$. Suppose $m\ne 0$. Since the Nakayama permutation of $A$ is cyclic, there is a number $n$ such that each indecomposable projective $A$-module is isomorphic to a direct summand of the terms of $\bigoplus_{1\leq s\leq n}\Lv_A^s(\X)$ in degrees $0$ and $m$. By Corollary \ref{basic},  $\bigoplus_{1\leq s\leq n}\Lv_A^s(\X)\simeq (\X)^{\oplus n}$ in $\Kb{A\mbox{\proj}}$. Since both $\bigoplus_{1\leq s\leq n}\Lv_A^s(\X)$ and $(\X)^{\oplus n}$ are radical complexes, it follows from \cite[(b), p.113]{HX} that $\bigoplus_{1\leq s\leq n}\Lv_A^s(\X)\simeq (\X)^{\oplus n}$ as complexes. Then each indecomposable projective $A$-module is isomorphic to a direct summand of $(X^{0})^{\oplus n}$ and $(X^{m})^{\oplus n}$. Thus $\Hom_{\Kb{A\mbox{\proj}}}((\X)^{\oplus n}, (\X)^{\oplus n}[m])\neq 0$. This contradicts to the fact that $(\X)^{\oplus n}$ is a tilting complex. Hence $m=0$ and $\X$ has only one nonzero term. $\square$

\bigskip
{\bf Acknowledgements}: The research work was supported partially by the National Natural Science Foundation of China (Grant 12031014 and 12226314). The authors thank Professor Wei Hu from Beijing Normal University, and Yiping Chen from Wuhan University for comments on the primary version of the manuscript. Also, the authors are grateful to Dr. Tiago Cruz from University of Stuttgart for pointing out an inaccurateness in the manuscript.

\end{document}